\documentclass[a4paper,10pt]{article}

\usepackage[latin1]{inputenc}
\usepackage{amsmath}
\usepackage{amsfonts}
\usepackage{amssymb}
\usepackage[english]{babel}
\usepackage{multirow}
\usepackage[dvips]{graphicx}
\usepackage{color}
\usepackage{indentfirst}
\usepackage{fancyhdr}

\newcommand{\Q}{\mathbb{Q}}
\newcommand{\qed}{\begin{flushright} \vspace{-0.5cm} $\blacksquare$ \end{flushright}}

\newtheorem{tm}{Theorem}[section]
\newtheorem{lemma}[tm]{Lemma}

\newtheorem{pp}[tm]{Proposition}

\newtheorem{cor}[tm]{Corollary}

\title{The Fermat-type equations $x^5 + y^5 = 2z^p$ or $3z^p$ solved through $\Q$-curves}
\author{Luis Dieulefait and Nuno Freitas}

\begin{document}

\maketitle

\begin{abstract} We solve the Diophantine equations $x^5 + y^5 = dz^p$ with $d=2, 3$ for a set of prime numbers of density $1/4$ , $1/2$, respectively. The method consists in relating a possible solution to another diophantine equation and solving the later by using $\mathbb{Q}$-curves and a generalized modular technique as in work of Ellenberg and Dieulefait-Jimenez along with some new techniques for eliminating newforms.
\end{abstract}

\section{Introduction}
Equations of the form
\begin{equation}
x^5 + y^5 = dz^p
\label{quintica}
\end{equation}
were first studied by Billerey in \cite{bil1} and then by Billerey-Dieulefait in \cite{BD} by using the classic modular approach via elliptic curves over $\mathbb{Q}$. Their work includes the solution for an infinite number of values for $d$ but the method there does not work in the cases $d=2$ or $d=3$. Here we will approach the problem for these values of $d$ making use of $\mathbb{Q}$-curves to apply a generalized modular technique as in the work of Dieulefait-Jimenez (see \cite{DJ}). To fulfill our strategy we also need to introduce and generalize techniques to eliminate newforms.\par
We call a triple $(a,b,c)$ such that $a^5 + b^5 = dc^p$ and $(a,b)=1$ a \textit{primitive} solution and we will say it is a \textit{trivial} solution if also $c=1$. In what follows we will show that for $d=2$ there is a subset of the prime numbers with density $1/4$ such that there are no non trivial primitive solutions to (\ref{quintica}) and that if $d=3$ the same is true but with a set of density $1/2$. In both cases for the rest of primes the non existence of solutions would follow immediately from a generalization of Ellenberg's work on the surjectivity of Galois representations attached to $\Q$-curves (see theorem 3.14 in \cite{ell2}) to real quadratic extensions.\par

Let $\gamma$ be an integer divisible only by primes $l \not\equiv 1$ (mod 5). In this paper we will prove the following theorems:

\begin{tm} For any $p > 13$ such that $p \equiv 1$ mod 4 and $p \equiv \pm 1$ mod 5, the equation $x^5 + y^5 = 2\gamma z^p$ has no non-trivial primitive solutions.
\label{tm1}
\end{tm}

\begin{tm} For any $p > 73$ such that $p \equiv 1$ mod 4, the equation $x^5 + y^5 = 3\gamma z^p$ has no non-trivial primitive solutions.
\label{tm2}
\end{tm}

Acknowledgements: We want to thank Nicolas Billerey for useful
discussions concerning the construction of the Frey
$\mathbb{Q}$-curve in section 3.

\section{A pair of diophantine equations}

From now on we will consider $a,b$ to be coprime integers, thus we will always be talking about primitive solutions. The first observation in order to solve equation (\ref{quintica}) is that we can relate a solution to (\ref{quintica}) with a solution of another Diophantine equation. Letting $\phi(a,b) = a^4 - a^3b + a^2b^2 - ab^3 + b^4$ we now recall some elementary results from \cite{bil1}:

\begin{lemma}The integers $a+b$ and $\phi(a,b)$ are coprime outside 5. Moreover, if 5 divides $a+b$ then $\upsilon_{5}(\phi(a,b))=1$
\label{coprimos}
\end{lemma}

\begin{lemma} Let $l \not\equiv 1$ (mod 5) be a prime number dividing $a^5 + b^5$. Then, $l$ divides $a+b$.
\label{dois}
\end{lemma}

From the fact that $a^5 + b^5 = (a+b)\phi(a,b)$ and the previous lemmas it follows easily that
\begin{cor} If a prime number $q \neq 5$ divides $\phi(a,b)$ then $q \equiv 1$ mod 5. Also, if $5 \nmid a+b$ then $\nu_5(\phi(a,b)) = 0$.
\label{coro}
\end{cor}

Let $(a,b,c)$ be a primitive solution to (\ref{quintica}) with $d$ fixed as in theorem \ref{tm1} or \ref{tm2}. Since $a^5 + b^5 = (a+b)\phi(a,b) = dc^p$, from the previous corollary it follows that $d \vert a+b$ and we must have either $\phi(a,b) = c_{0}^{p}$ or $\phi(a,b) = 5c_{0}^{p}$ if $5 \nmid a+b$ or $5 \vert a+b$, respectively, where $c_{0} \vert c$ is an integer coprime with $d$ and 5. Hence we will prove that (\ref{quintica}) has no primitive solutions if we show that there are no primitive solutions to both equations
\begin{equation}
\phi(a,b) = c^p,
\label{caso1}
\end{equation}
and
\begin{equation}
\phi(a,b) = 5c^p
\label{caso5}
\end{equation}
such that $d \vert a+b$ and $d$, 5 coprime with $c$. Also observe that by the corollary above all the primes dividing $c$ are congruent with 1 modulo 5. Then, in particular, $c$ is odd. Although theorem \ref{tm1} and theorem \ref{tm2} are true for an infinite set of values for $d$ in what follows we will suppose that $d=2$ or $d=3$, because this already implies $2 \vert a+b$ or $3 \vert a+b$ and that is all we need to prove theorem \ref{tm1} and theorem \ref{tm2}, respectively.\\

\section{The $\mathbb{Q}$-curve}

To apply the modular approach we need to find an appropriate Frey $\Q$-curve.
Consider the polynomial
$$\phi(x,y) = x^4 - x^3y + x^2y^2 - xy^3 + y^4$$
and its factorization over $\mathbb{Q}(\sqrt{5})$
$$\phi(x,y) = \phi_{1}(x,y) \phi_{2}(x,y),$$ where
\begin{eqnarray*}
\phi_{1}(x,y) &=& x^2 + \omega xy + y^2, \\
\phi_{2}(x,y) &=& x^2 + \bar{\omega} xy + y^2, \\
\omega &=& \frac{-1+\sqrt{5}}{2}, \\
\bar{\omega}&=&\frac{-1-\sqrt{5}}{2}.
\end{eqnarray*}

\begin{pp} If $(a,b)$ is a pair of coprime integers then $\phi_1(a,b)$ and $\phi_2(a,b)$ are coprime outside 5.
 \label{ppot}
\end{pp}
\textbf{Proof:} Suppose that $l$ is a prime in $\Q(\sqrt{5})$ dividing both $\phi_i(a,b)$, then also divides $\phi_1 - \phi_2 = \sqrt{5}ab$. If $l$ divide $ab$, then we can suppose that $l$ divides $a$, but dividing $a$ and $\phi_1(a,b)$ implies that divides $b$ which is a contradiction since $a,b$ are coprime. Thus we conclude that $l$ is above 5.\qed

Let $(a,b)$ be a non-trivial primitive solution to (\ref{caso1}) or (\ref{caso5}), we define an elliptic curve over $\mathbb{Q}(\sqrt{5})$
$$E_{(a,b)}: y^2 = x^3 + 2(a+b)x^2 - \bar{\omega}\phi_{1}(a,b)x.$$ To ease notation
we will omit the pair $(a,b)$ when writing $E_{(a,b)}$, $\phi(a,b)$, $\phi_{1}(a,b)$, $\phi_{2}(a,b)$. Observing that
$(a+b)^2 = -\bar{\omega}\phi_1 -\omega\phi_2$ we compute
\begin{equation}
\Delta(E) = 2^6\bar{\omega}\phi\phi_{1}.
\label{disc}
\end{equation}
Consider the Galois conjugated curve for the non trivial element in
$\mbox{Gal}(\mathbb{Q}(\sqrt{5})/\mathbb{Q})$
$${^{\sigma}}E: y^2 = x^3 + 2(a+b)x^2-\omega\phi_{2}x,$$
and the 2-isogeny $\mu : {^{\sigma}}E \rightarrow E$ given by
$$(x,y) \mapsto (-\frac{y^2}{2x^2},\frac{\sqrt{-2}}{4}\frac{y}{x^2}(\omega\phi_2 + x^2)).$$
with dual isogeny $\hat{\mu} : E \rightarrow {^{\sigma}}E$ given by
$$(x,y) \mapsto (-\frac{y^2}{2x^2},- \frac{\sqrt{-2}}{4}\frac{y}{x^2}(\bar{\omega}\phi_1 + x^2)),$$
showing that $E$ is a $\Q$-curve. In order to apply the modular approach to our equation, we need to find a
twist of $E$, $E_{\gamma}$ such that its Weil restriction decomposes as a product of abelian varieties of $GL_2$-type. That is the content of the following theorem, whose proof we postpone until section \ref{prova}.
\begin{tm}Let $K=\Q(\theta)$ where $\theta = \sqrt{\frac{1}{2}(5 + \sqrt{5})}$ is a root of the
polynomial $x^4 - 5x^2 + 5$, and put $\gamma = 2\theta^2 - \theta - 5$. Then, the twisted curve
$$E_{\gamma} : y^2 = x^3 + 2\gamma(a+b)x^2 -\gamma^2\bar{\omega}\phi_{1}(a,b)x$$
is completely defined over $K$ and its Weil restriction decomposes as the product of two non-isogenous abelian surfaces of $GL_2$-type each of them with endomorphism algebra over $\Q$ isomorphic to $\Q(i)$.
\label{twist}
\end{tm}

\subsection{The conductor of $E_{\gamma}$}

In this subsection we list the results about the conductor of $E_{\gamma}$ which are needed to compute the conductor of its Weil restriction and consequently the level of the associated modular forms. All the computations that are required to obtain these results were made following Tate's algorithm as in section 9 of chapter IV of \cite{sil2}.\par
Let $K=\Q(\theta)$ and denote by $\mathfrak{P}_{2}$ and $\mathfrak{P}_{5}$ the only primes above 2 and 5, respectively. Also, denote the conductor of $E_{\gamma}$ by $N_{E_{\gamma}}$ and let $\mbox{rad}(c) = \prod\limits_{p \vert c} p$ be the product of the prime factors of $c$. We have different results on $N_{E_{\gamma}}$ for the case $d=2$ or $d=3$. First consider $d=2$ and observe that it implies $a+b$ to be even, hence the following holds.
\begin{pp} If $(a,b,c)$ is a primitive solution of equation (\ref{caso1}) coming from $d=2$ then the curve $E_{\gamma}$ has conductor:
$$N_{E_{\gamma}} = \left\{ \begin{array}{ll}
\mathfrak{P}_{2}^8\mathfrak{P}_{5}^2({\mbox{rad}(c)}) & \mbox{ if } \quad 2 \Vert a+b, \vspace{0.1cm} \\
\mathfrak{P}_{5}^{2}(\mbox{rad}(c)) & \mbox{ if } \quad 4 \Vert a+b, \vspace{0.1cm} \\
\mathfrak{P}_{2}^4\mathfrak{P}_{5}^2({\mbox{rad}(c)}) \quad & \mbox{ if } \quad 8 \vert a+b.
\end{array} \right.$$
\end{pp}

\begin{pp} If $(a,b,c)$ is a primitive solution of equation (\ref{caso5}) coming from $d=2$ then the curve $E_{\gamma}$ has conductor:
$$N_{E_{\gamma}} = \left\{ \begin{array}{ll}
\mathfrak{P}_{2}^8{(\mbox{rad}(c))} & \mbox{ if } \quad 2 \Vert a+b, \vspace{0.1cm} \\
(\mbox{rad}(c)) & \mbox{ if } \quad 4 \Vert a+b \vspace{0.1cm} \\
\mathfrak{P}_{2}^4{(\mbox{rad}(c))} & \mbox{ if } \quad 8 \vert a+b.
\end{array} \right.$$
\end{pp}

In the case $d=3$ the solutions to (\ref{quintica}) may verify $a+b$ odd. This results in an extra case of the conductor at $\mathfrak{P}_{2}$.

\begin{pp} If $(a,b,c)$ is a primitive solution of equation (\ref{caso1}) coming from $d=3$ with $a+b$ odd then the curve $E_{\gamma}$ has conductor:
$$N_{E_{\gamma}} =  \mathfrak{P}_{2}^8\mathfrak{P}_{5}^2{(\mbox{rad}(c))} \quad \mbox{ or } \quad N_{E_{\gamma}} = \mathfrak{P}_{2}^{6}\mathfrak{P}_{5}^{2}(\mbox{rad}(c)).$$
\end{pp}

\begin{pp} If $(a,b,c)$ is a primitive solution of equation (\ref{caso5}) coming from $d=3$ with $a+b$ odd then the curve $E_{\gamma}$ has conductor:
$$N_{E_{\gamma}} =  \mathfrak{P}_{2}^8{(\mbox{rad}(c))} \quad \mbox{ or } \quad N_{E_{\gamma}} = \mathfrak{P}_{2}^{6}(\mbox{rad}(c)).$$
\end{pp}

Let $E_{\gamma,2}$ be the twist of $E_{\gamma}$ by 2 and denote its conductor by $N_{E_{\gamma,2}}$. This conductor will later be used to reach a contradiction when analyzing representations coming from newforms. Actually, we will only need it to eliminate newforms of level 1600 for the case $d=2$, in which case we can suppose that $2 \Vert a+b$.

\begin{pp} Suppose that $(a,b,c)$ is a solution to $(\ref{caso1})$ or $(\ref{caso5})$ such that $2 \Vert a+b$. Then the conductor at $\mathfrak{P}_{2}$ of $E_{\gamma,2}$ is $\mathfrak{P}_{2}^{0}$ or $\mathfrak{P}_{2}^4$.
\label{cond22}
\end{pp}
\label{cond}

\section{Modularity of $E_{\gamma}$}

Here we go into determining the precise Serre's parameters: weight, level and character of a residual representation of a $GL_2$-type abelian variety attached to our Frey-curves.\\

Let $B = \mbox{Res}_{K/\Q}(E_{\gamma}/K)$, where $K=\Q(\theta)$ is the cyclic Galois extension as before, and denote its conductor by $N_{B}$. To compute this conductor we will use a formula of Milne (see \cite{mil1}), which in our case tells us that the conductor satisfies
$$N_{B} = \mbox{Nm}_{K/\Q}(N_{E_{\gamma}})\mbox{Disc}(K/\Q)^2,$$ where Disc and Nm denote the discriminant and norm of $K/\Q$, respectively. Since $c$ is odd and $\mbox{Disc}(K/\Q) = 2^{4}5^{3}$, then the primes dividing $c$ do not ramify in $K$. Being $K/\Q$ of degree 4 we have that a prime $p$ dividing $c$ is inert, splits completely or is the product of two primes in $K$ with residual degree 2 thus $\mbox{Nm}((c))=c^4$. Also, $\mbox{Nm}(\mathfrak{P}_{5})=5$ and $\mbox{Nm}(\mathfrak{P}_{2})=4$. By applying the above formula we obtain the following results:
\begin{pp} If $(a,b,c)$ is a primitive solution of equation (\ref{caso1}) then $B$ has conductor:
$$N_{B} = \left\{ \begin{array}{ll}
2^{24}5^8\mbox{rad}(c)^4 \mbox{ or  } 2^{20}5^8\mbox{rad}(c)^4 & \mbox{ if } \quad 2 \nmid a+b, \vspace{0.1cm} \\
2^{24}5^8\mbox{rad}(c)^4 & \mbox{ if } \quad 2 \Vert a+b, \vspace{0.1cm} \\
2^{8}5^8\mbox{rad}(c)^4 & \mbox{ if } \quad 4 \Vert a+b, \vspace{0.1cm} \\
2^{16}5^8\mbox{rad}(c)^4 \quad & \mbox{ if } \quad 8 \vert a+b.
\end{array} \right.$$
\end{pp}

\begin{pp} If $(a,b,c)$ is a primitive solution of equation (\ref{caso5}) then $B$ has conductor:
$$N_{B} = \left\{ \begin{array}{ll}
2^{24}5^6\mbox{rad}(c)^4 \mbox{ or  } 2^{20}5^6\mbox{rad}(c)^4 & \mbox{ if } \quad 2 \nmid a+b, \vspace{0.1cm} \\
2^{24}5^6\mbox{rad}(c)^4 & \mbox{ if } \quad 2 \Vert a+b, \vspace{0.1cm} \\
2^{8}5^6\mbox{rad}(c)^4 & \mbox{ if } \quad 4 \Vert a+b, \vspace{0.1cm} \\
2^{16}5^6\mbox{rad}(c)^4 \quad & \mbox{ if } \quad 8 \vert a+b.
\end{array} \right.$$
\end{pp}

We know from theorem \ref{twist} that $B \simeq S_{1} \times S_{2}$ where $S_{i}$ are two non $\Q$-isogenous abelian surfaces of $GL_{2}$-type with $\Q$-endomorphism algebras equal to $\Q(i)$. So the conductor of $B$ satisfies $N_{B}=N_{S_{1}}N_{S_{2}}$.\par

For a prime $l$ and each $S_{i}$ the action on the Tate module $T_{l}S_i$ induces a 4-dimensional $l$-adic representation of $G_{\Q}$ that decomposes into two 2-dimensional $\lambda$-adic representations $\rho_{S_i,\lambda}$ and $\rho_{S_i,\lambda}^{\sigma}$, where $\lambda$ is a prime of $\Q(i)$ above $l$. Then we have four 2-dimensional representations of $G_{\Q}$ extending the $l$-adic representation $\rho_{E_{\gamma},l}$ of $\mbox{Gal}(\bar{\Q}/K)$ induced by the action on the Tate module of $E_{\gamma}$. Since extensions of absolutely irreducible representations are unique up to twists we have the following relations between them:
$$\rho_{S_1,\lambda} \otimes \epsilon = \rho_{S_1,\lambda}^{\sigma}, \quad \quad \rho_{S_1,\lambda} \otimes \epsilon^2 = \rho_{S_2,\lambda} \quad \mbox{ and } \quad \rho_{S_1,\lambda} \otimes \epsilon^3 = \rho_{S_2,\lambda}^{\sigma},$$
where $\epsilon$ is the character of $K$ (see section \ref{prova}, formula (\ref{caracter})) and $\epsilon^2$ is the character of $\Q(\sqrt{5})$. It is known that the conductors of $\rho_{S_i,\lambda}$ and $\rho_{S_i,\lambda}^{\sigma}$ are equal and that their product is equal to the conductor of $S_i$, which means that every prime in the conductor of $S_i$ appears to an even power. From the second relation and the fact that $\epsilon^2$ has conductor 5 we see that the only possible difference in the conductors of $\rho_{S_1,\lambda}$ and $\rho_{S_2,\lambda}$ may occur at 5. Furthermore, the conductor at 5 of $\rho_{S_1,\lambda} \otimes \epsilon^2$ is equal to the least common multiple between the conductor at 5 of $\rho_{S_1,\lambda}$ and $\mbox{cond}(\epsilon^2)^2 = 5^2$, except if these two numbers coincide in which case the twist may have a smaller conductor. Using these facts together with a case checking allow us to determine all the possibilities for the conductors of the four 2-dimensional representations (see table \ref{conductors}).\\
\begin{table}[ht]
\begin{center}
\begin{tabular}{|c|c|c|c|c|c|}
\hline
\mbox{Equation}& $\nu_{2}(a+b)$& $\rho_{S_1,\lambda}$& $\rho_{S_1,\lambda}^{\sigma}$& $\rho_{S_2,\lambda}$&$\rho_{S_2,\lambda}^{\sigma}$\\ \hline
(\ref{caso1})  & 0              & $2^{6}5^{2}c_0$       & $2^{6}5^{2}c_0$             & $2^{6}5^{2}c_0$       & $2^{6}5^{2}c_0$\\ \hline
(\ref{caso1})  & $\geq 3$           & $2^{6}5^{2}c_0$       & $2^{6}5^{2}c_0$             & $2^{6}5^{2}c_0$       & $2^{6}5^{2}c_0$\\ \hline
(\ref{caso1})  & 2              & $2^{4}5^{2}c_0$       & $2^{4}5^{2}c_0$         & $2^{4}5^{2}c_0$       & $2^{4}5^{2}c_0$\\ \hline
(\ref{caso1})  & 1              & $2^{2}5^{2}c_0$       & $2^{4}5^{2}c_0$             & $2^{4}5^{2}c_0$       & $2^{4}5^{2}c_0$ \\ \hline
(\ref{caso5})  & 0              & $2^{5}5^{2}c_0$       & $2^{5}5^{2}c_0$             & $2^{5}5c_0$           & $2^{5}c_0$\\ \hline
(\ref{caso5})  & $\geq 3$           & $2^{6}5^{2}c_0$       & $2^{6}5^{2}c_0$             & $2^{6}5c_0$           & $2^{6}5c_0$  \\ \hline
(\ref{caso5})  & 2                  & $2^{4}5^{2}c_0$       & $2^{4}5^{2}c_0$             & $2^{4}5c_0$           & $2^{4}5c_0$  \\ \hline
(\ref{caso5})  & 1              & $2^{2}5^{2}c_0$       & $2^{2}5^{2}c_0$             & $2^{2}5c_0$           & $2^{2}5c_0$  \\ \hline
\end{tabular}
\caption{Values of conductors, where $c_0 = \mbox{rad}(c)$}
\label{conductors}
\end{center}
\end{table}

Now pick a prime $\lambda$ in $\Q(i)$ above $p$ and let $\bar{\rho} := \bar{\rho}_{S_1,\lambda}$ be the residual representation of $\rho_{S_1,\lambda}$. Recall that the Frey-Hellegouarch argument (used in the proof of FLT) states that a prime of semistable reduction of an elliptic curve $C$ which appear in the discriminant of $C$ to a $p$-power will not ramify in the mod $p$ representation induced by the $p$-torsion points. From formula (\ref{disc}), proposition \ref{ppot} and the fact that the primes in $K$ dividing $c$ are of semistable reduction for $E_{\gamma}$ we can apply the Frey-Hellegouarch argument to conclude that the restriction $\bar{\rho}_{|K}$ of $\bar{\rho}$ to $\mbox{Gal}(\bar{\Q}/K)$, which coincides with $\bar{\rho}_{E_{\gamma},p}$, will not ramify at primes dividing $c$. Since $K$ only ramifies at 2 and 5 then we see that $\bar{\rho}$ can not ramify outside 2 and 5. On the other hand it is well known that in the presence of wild ramification the conductor does not decrease when reducing mod $p$, so the possible conductors of $\bar{\rho}$ are exactly the values in the third column of table \ref{conductors} without the factor $c_0$. Hence we have determined Serre's level $N(\bar{\rho})$ needed to apply Serre strong conjecture (which is now a theorem). The following two propositions tell us the Serre weight and character.

\begin{pp}The representation $\bar{\rho}$ has character $\bar{\epsilon}$ (conjugate of $\epsilon$).
\end{pp}
\textbf{Proof:} By Theorem 5.12 of \cite{pyle} we know that the character associated to $\rho_{S_1,\lambda}$ is equal to $\epsilon^{-1}$ where $\epsilon$ is the splitting character defined  by formula (\ref{caracter}) in section \ref{prova}. Furthermore, since $\epsilon$ has order 4, for any prime $p$ distinct from 2 the representation $\bar{\rho}$ has the same character of the non residual one. Then $\bar{\rho}$ is modular with character $\epsilon^{-1} = \bar{\epsilon}$.\qed

\begin{pp}The Serre weight of the representation $\bar{\rho}$ is $k(\bar{\rho})=2$.
\end{pp}
\textbf{Proof:} We divide in two cases. If $p \nmid c$ then $S_1$ has good reduction at $p$ and $\bar{\rho}$ comes from an abelian variety with good reduction at $p$ hence $k=2$. If $p \vert c$, the fact that $p\vert \nu_{\mathfrak{P}}(\Delta)$ for any $\mathfrak{P}$ above $p$ implies as in \cite{ell2} that $\bar{\rho}$ is finite and so $k=2$.\qed

Finally, to apply the Serre's conjecture we still need irreducibility.

\begin{pp}The representation $\bar{\rho}$ is absolutely irreducible.
\end{pp}
\textbf{Proof:} If $(a,b,c)$ is a non-trivial solution then there is a prime of characteristic greater that 3 of semistable reduction, hence by proposition 3.2 in \cite{ell2} we conclude that $\bar{\rho}$ irreducible for $p > 13$.\qed

Now from Serre's strong conjecture we know that there must be a newform $f$ in the space $\mathcal{S}_{2}(M,\bar{\epsilon})$ with $M=1600$, 800, 400 or 100 and $\mathfrak{P} \vert p$ in $\Q_{f}$ such that the mod $\mathfrak{P}$ residual representation attached to $f$ (denoted by $\bar{\rho}_{f,p}$) is isomorphic to our residual representation $\bar{\rho}$.

\section{Eliminating Newforms}

Using the SAGE software we compute the newforms in the spaces determined in the previous section. We need to show that our representation $\bar{\rho}$ can not come from any of the computed newforms in order to reach a contradiction.\par
Note that equations $\phi(a,b) = \pm 1$ and $\phi(a,b) = \pm 5$ have solutions $(\pm 1,0)$, $(0,\pm 1)$, $(1,1)$, $(-1,-1)$ and $(1,-1)$, $(-1,1)$, respectively. This means that for these pairs $(a,b)$ the Frey $\Q$-curves do exist and so if their corresponding newforms have level 100, 400, 800 or 1600 a priori we may not be able to eliminate those forms. Recall that when treating the case $d=2$ or $d=3$ we have $a+b$ even or $3 \mid a+b$, respectively. As we will see this will turn out to be a key information. Now observe that $(\pm 1,0)$, $(0, \pm 1)$ will not be a problem for equation (\ref{quintica}) with $d=2$, because $a+b$ is odd; also $E_{(-1,1)}$ and $E_{(1,-1)}$ correspond to level 320. Actually, only $E_{(1,1)}$ and $E_{(-1,-1)}$ correspond to newforms of level 1600 which we will need to eliminate. We will be able to do this because these curves correspond to newforms with complex multiplication.\\

\subsection{The case $d=2$}

Since $d=2$ we have $a+b$ even thus $800 = 2^5 5^2 $ is not a possible level. To eliminate the newforms in the other levels we separate them into three sets and then we will apply a different strategy for each of the sets. Observe that $\Q(i) = \Q(\epsilon) \subseteq \Q_{f}$ and let S1, S2 and S3 be as follows:
\begin{itemize}
\item[\textbf{S1:}] Newforms with CM (Complex Multiplication),
\item[\textbf{S2:}] Newforms without CM and field of coefficients strictly containing $\Q(i)$,
\item[\textbf{S3:}] Newforms without CM and field of coefficients $\Q(i)$
\end{itemize}
We now derive a contradiction for each set.\\

\textbf{S1:} Modulo Galois conjugation there are 8 newforms with complex multiplication. Half of them with CM by $\Q(i)$ and the other half by $\Q(\sqrt{-5})$. If there is a non-trivial solution $(a,b,c)$ to equation (\ref{caso1}) or (\ref{caso5}) then there exists a prime not dividing 6 of semistable reduction for $E$. Hence if $p > 13$ by prop. 3.4 in \cite{ell2} the images of $\bar{\rho}$ will not lie in the normalizer of a split Cartan subgroup. Then if $p$ is a square in both $\Q(i)$ and $\Q(\sqrt{-5})$, i.e. $p \equiv 1$ mod 4 and $p \equiv \pm 1$ mod 5 we can not have $\bar{\rho} \equiv \bar{\rho}_{f,p}$ for $f$ in S1. This is because for a newform with CM it is known that for $p$'s that are squares on the field of complex multiplication the image of the attached representation  will be in a normalizer of a split Cartan subgroup.\\

\textbf{S2:} There are 12 newforms in this group (modulo conjugation). From propositions 3.4 and 4.3 in \cite{rib1} we know that for the primes $q$ of good reduction for $S_1$ the values $a_{q}$ satisfy $a_q = \bar{a}_q\bar{\epsilon}(q)$, i.e. an inner twist coming from the nebentypus. Here $\epsilon$ is the character of order 4 defined in section \ref{prova}; in particular, the inner twist implies that
$$a_q = \left\{ \begin{array}{ll}
t & \mbox{ if } q \equiv 1 \mbox{ or } 19 \mbox{ mod } 20 \\
it & \mbox{ if } q \equiv 9 \mbox{ or } 11 \mbox{ mod } 20  \\
t - it & \mbox{ if } q \equiv 3 \mbox{ or } 17 \mbox{ mod } 20 \\
t + it & \mbox{ if } q \equiv 7 \mbox{ or } 13 \mbox{ mod } 20
\end{array} \right.$$
where $t$ is an integer. From corollary \ref{coro} it follows that 3 is a prime of good reduction of $S_1$. Hence, $a_3$ must be of the form $t-ti$ and from the Weil bound $|a_q| \leq 2\sqrt{q}$ follows that $|t|\leq 2$. If $f$ is one of these 12 newforms then the congruence
\begin{equation}
a_3 \equiv c_3(f) \mbox{ (mod } \mathfrak{P})
\label{cong}
\end{equation}
must hold, for some newform $f = q + \sum_{n=2} {c_n}q^n$ in S2 and a prime $\mathfrak{P}$ in $\bar{\Q}$ above $p$. Now for each newform in S2 we use the coefficient $c_3$ to derive a contradiction, because none of them has $c_3$ of the form $t-ti$. As an example, there is a newform in S2 of level 400 with $c_3$ having minimal polynomial $x^2 +10i$; thus if the congruence (\ref{cong}) holds we must have $c_3 \equiv t-it \mbox{ (mod } \mathfrak{P})$ with $t=0,\pm 1,\pm 2$. Taking fourth powers we get $100 \equiv 4t^4 \mbox{ (mod } \mathfrak{P})$ which means $25 - t^4 \equiv 0 \mbox{ (mod } \mathfrak{P})$ and finally, substituting for the possible values of $t$ we reach a contradiction if $p >5$. A similar argument works for every newform in S2 and we conclude that if $p > 7$ the newform predicted by Serre's conjecture can not be in S2.\\

\textbf{S3:} There are 10 newforms in this group all with level 1600. Recall that $\bar{\rho} := \bar{\rho}_{S_1,\lambda}$  and suppose that $\bar{\rho} \equiv \bar{\rho}_{f, p}$ for some $f$ in S3. A priori all newforms in S3 are inconvenient, in the sense that their Fourier coefficients ${a_{q}}'s$ behave exactly as those of our surface $S_1$. That is, the ${a_q}'s$ lie in $\Q(i)$, they respect the rule $a_q = \bar{a}_q\bar{\epsilon}(p)$ (by construction) and they are even (this is true for our surface because $E_{\gamma}$ has a 2-torsion point). To deal with this problem we will twist each newform in S3 by the character of $\Q(\sqrt{2})$, which we denote by $\chi$. Note that $\mbox{cond}(\chi)^2=8^2=2^6$ and the power of 2 in 1600 is also $2^6$ so as mentioned before we are in a situation where the power of 2 in the level of the twist can decrease. Indeed, using SAGE, we pick $f$ in S3, we consider $f \otimes \chi$ and compare the coefficients of $f \otimes \chi$ up to the Sturm bound with those of the newforms with level dividing 1600 to find that $f \otimes \chi$ is a newform of level 800 for all $f$ in S3.\par
On the other hand, let $E_{\gamma,2}$ be as in section \ref{cond} and $\rho_{E_{\gamma,2},p}$ be the representation coming from the action of $G_{\Q}$ on the Tate module $T_{p}E_{\gamma,2}$. Note that
$$(\rho_{S_1,\lambda} \otimes \chi)_{|K} = (\rho_{S_1,\lambda})_{|K} \otimes \chi_{|K} = \rho_{E_{\gamma},p} \otimes \chi_{|K},$$ that is $\rho_{S_1,\lambda} \otimes \chi$ extends $\rho_{E_{\gamma},p} \otimes \chi_{|K}$ and this one is precisely $\rho_{E_{\gamma,2},p}$. The same is true for the other three representations coming from the Weil restriction $B \simeq S_1 \times S_2$. Moreover, $\rho_{B,p} = \mbox{Ind}_{G_{K}}^{G_{\Q}}{\rho_{E_{\gamma},p}}$ and we have
$$\rho_{B,p} \otimes \chi = (\mbox{Ind}_{G_{K}}^{G_{\Q}}{\rho_{E_{\gamma},p}}) \otimes \chi = \mbox{Ind}_{G_{K}}^{G_{\Q}}{(\rho_{E_{\gamma},p} \otimes \chi_{|K}}) = \mbox{Ind}_{G_{K}}^{G_{\Q}}{\rho_{E_{\gamma,2},p}},$$
and this means that $\rho_{B,p} \otimes \chi$ is the representation coming from the action of $G_{\Q}$ on the Tate module of $\mbox{Res}_{K/\Q}(E_{\gamma,2}/K)$. Let $\rho_{1}$ denote the 2-dimensional factor $\rho_{S_1,\lambda} \otimes \chi$ of $\rho_{B,p} \otimes \chi$. From proposition \ref{cond22}, Milne's formula and an analysis identical to that used to compute table \ref{conductors} we conclude that $N(\bar{\rho}_1) = 400$ or 100. $\bar{\rho}_1$ is also irreducible and has character and Serre's weight equal to those of $\bar{\rho}$, because the same arguments hold. We now apply Serre's strong conjecture to $\bar{\rho}_{1}$ and we conclude that there must be a newform $g$ in level 400 or 100 such that $\bar{\rho}_{1} \equiv \bar{\rho}_{g,p}$ (mod $\mathfrak{P}$). But at the same time we also have $$\bar{\rho}_{1} = \overline{\rho_{S_1,\lambda} \otimes \chi} = \bar{\rho} \otimes \chi \equiv \bar{\rho}_{f,p} \otimes \chi \equiv \bar{\rho}_{f \otimes \chi,p} = \bar{\rho}_{f',p},$$
where all congruences are mod $\mathfrak{P}$ and $f'$ is of level 800 by the previous paragraph. Hence the congruence $$\bar{\rho}_{g,p} \equiv \bar{\rho}_{f',p} \mbox{ (mod } \mathfrak{P})$$ must hold between a newform of level 800 and another of level 400 or 100, but Carayol has shown (see \cite{car}) that this kind of level lowering can not happen. This ends the proof of theorem \ref{tm1}.

\subsection{The case $d=3$}

In this case $a+b$ may be odd and as it was seen before we now need to consider the extra level 800, but in our favor we always have $3 \vert a+b$. We start by considering the case of $a+b$ odd which leaves us with newforms of level 1600 or 800.\par
In level 800 there are 4 newforms of type S2 and 10 of type S3 and none of type S1. The same type of argument as for $d=2$ will eliminate the forms of type S2, but this time we have a contradiction only for $p > 73$, because there are newforms such that $a_3$ has minimal polynomial $t^2 \pm (2 - 2i)t + i$. To eliminate a newform $f$ of type S3 (level 800), let $\mathfrak{P}_3$ be the prime of $K$ above 3, and  suppose that $\bar{\rho} \equiv \bar{\rho}_{f,p}$ (mod $\mathfrak{P}$). In particular, $\bar{\rho}|_{G_{K}} \equiv \bar{\rho}_{f,p}|_{G_{K}}$ and so the values $a_{\mathfrak{P}_3}(E_{\gamma})$ and $a_{\mathfrak{P}_3}(f)$ must also be congruent mod $p$. Now, on one hand it is easy to check by direct computation that the hypothesis $3 \vert a+b$ implies that $E_{\gamma}$ mod $\mathfrak{P}_{3}$ is the same for any pair $(a,b)$ and that $E_{\gamma}$ has $a_{\mathfrak{P}_{3}}(E_{\gamma}) = -18$ for all $a,b$. On the other hand looking for the $a_3$ coefficients of the newforms of type S3 we find four possibilities $\pm(2i - 2)$ and $\pm(i-1)$. It is known that if $\alpha$, $\beta$ are the roots of the polynomial $x^2 - a_{3}x + \bar{\epsilon}(3)3$ (the characteristic polynomial of $\rho_{f,p}(\mbox{Frob}_{3})$), then $a_{\mathfrak{P}_{3}}(f) = \alpha^4 + \beta^4$. Since $\bar{\epsilon}(3) = -i$, by substituting for each of the four values of $a_3$ we find out that $a_{\mathfrak{P}_{3}}(f) = 14$ or 2. Hence $\bar{\rho} \equiv \bar{\rho}_{f,p}$ is not possible if $p > 3$.\par
In level 1600 the forms of type S1 and S2 can be eliminated exactly as for $d=2$, but for type S1 we get a better result because newforms with complex multiplication by $\Q(\sqrt{-5})$ verify $a_3 = \pm(i -1)$ which we have just seen to be in contradiction with $3 \vert a+b$, hence to complete the contradiction we only need to suppose that $p$ is a square in $\Q(i)$, that is, $p \equiv 1$ mod 4. For those $f$ of type S3 we first twist them by the Dirichlet character of $\Q(\sqrt{2})$ and we already know that $f \otimes \chi$ is a newform of level 800. On the other hand we twist the Frey curve $E_{\gamma}(a,b)$ by the same character. If the conductor at 2 of the twisted curve is not $2^6$ we have a contradiction via Carayol as in the case $d=2$; if it is $2^6$, then since $E_{\gamma} \otimes \chi$ mod $\mathfrak{P}_3$ is equal to $E_{\gamma}(a,b)$ mod $\mathfrak{P}_3$ we are in the same situation as above with forms of level 800 and the argument used there with the values of $a_{\mathfrak{P}_3}(E_{\gamma})$ and $a_{\mathfrak{P}_3}(f)$ gives us the desired contradiction.\par
For the remaining cases, those with $a+b$ even, we eliminate newforms of type S2 and S3 exactly with the same arguments used when $d=2$, and again for those of type S1 we only need to suppose that is $p \equiv 1$ mod 4 to get a contradiction, because the newforms with complex multiplication by $\Q(\sqrt{-5})$ verify $a_3 = \pm (i-1)$. This ends the proof of theorem \ref{tm2}.

\section{Finding $E_{\gamma}$}
\label{prova}

In this section we prove theorem \ref{twist} by extensively using the work of Quer, see \cite{quer1} and \cite{quer2} for the notation and theorems. First note that the theory of Quer always require the curves not to have complex multiplication, hence we need
\begin{tm} Let $(a,b,c)$ be a primitive non trivial solution of equation (\ref{caso1}) or (\ref{caso5}). Then the curve $E = E_{(a,b)}$ has no complex multiplication.
\end{tm}
\textbf{Proof:} This theorem follows immediately from the fact that the primes dividing $c$ are of semistable reduction of $E$.\qed
The minimal field of definition of $E$ is $K_d = \Q(\sqrt{5})$ and $\{a_1\} = \{5\}$, $\{d_1\} = \{2\}$ is a dual base respect to the corresponding degree map. Since $(5,\pm 1) \neq 1$ from proposition 4.3b of \cite{quer1} we know that $K_d$ is not a splitting field. Let $L=\Q(\sqrt{5},\sqrt{-2})$ be a biquadratic field of complete definition of $E$ and let $\sigma , \tau \in Gal(L/\Q)$ be such that  $$\left\{ \begin{array}{lll}
                     \sigma(\sqrt{5}) = \sqrt{5}  & \mbox{ and } & \sigma(\sqrt{-2}) = -\sqrt{-2} \\
                     \tau(\sqrt{5}) = -\sqrt{5}  & \mbox{ and } & \tau(\sqrt{-2}) = \sqrt{-2}\
                  \end{array} \right.$$
From the explicit expression of the isogeny we compute the 2-cocycle $$c_{L}(g,h) = \phi_g{^{g}}\phi_h\phi_{gh}^{-1}$$ obtaining
\begin{table}[h]
\begin{center}
\begin{tabular}{cc||c|c|c|c|}
\cline{3-6}
& & \multicolumn{4}{|c|}{$h$} \\ \cline{3-6}
& & 1 & $\sigma$ & $\tau$ & $\sigma\tau$ \\ \hline \hline
\multicolumn{1}{|c|}{\multirow{4}{*}{$g$}} &
1 & 1 & 1 & 1 & 1      \\ \cline{2-6}
\multicolumn{1}{|c|}{}                        &
$\sigma$ & 1 & 1 & -1 & -1      \\ \cline{2-6}
\multicolumn{1}{|c|}{} &
$\tau$ & 1 & 1 & -2 & -2  \\ \cline{2-6}
\multicolumn{1}{|c|}{}                        &
$\sigma\tau$ & 1 & 1 & 2 & 2  \\ \hline
\end{tabular}
\caption{Values of $c_{L}$}
\label{cociclo}
\end{center}
\end{table}

From the table we see that we need an appropriate splitting character in order to apply theorem 4.2 and proposition 5.2 in \cite{quer1}. First by theorem 3.1 in \cite{quer1} one sees that $\xi(C)_{\pm} = (5,2)$ and using local information it can be checked that the character $\epsilon: (\mathbb{Z}/20\mathbb{Z})^{*} \rightarrow \Q(\zeta_{4})$ of order 4 and conductor 20 given by
\begin{equation}
\epsilon=\epsilon_{2}\epsilon_{5},
\label{caracter}
\end{equation}
where $\epsilon_2$ is the quadratic character of conductor 4 and $\epsilon_{5}$ is the character of order 4 and conductor 5 given by $\epsilon_{5}(2)=\zeta_4=i$ will be an appropriate splitting character. Its fixed field is $K_{\epsilon}=\Q(\theta)$, where $\theta = \sqrt{\frac{1}{2}(5 + \sqrt{5})}$. Then it follows from theorem 4.2 and proposition 5.2 in \cite{quer1} that there is a representative in the isogeny class of $E$ which is completely defined over $K_{\beta} = K_{\epsilon}K_{d}=\Q(\theta)\Q(\sqrt{5})=\Q(\theta)$ and has trivial Schur class. We denote it by $E_{\gamma}$.\par
Now we will explicitly solve an embedding problem to find $\gamma$. As explained in section 3 of \cite{quer2}, let
\begin{equation}
\label{soluvel}
\xi_{K_{\beta}L} = (\mbox{Inf}_{Gal(K_{\beta}/\Q)}^{Gal(K_{\beta}L/\Q)}c_{K_{\beta}}(E_{\gamma}))_{\pm}(\mbox{Inf}_{Gal(L/\Q)}^{Gal(K_{\beta}L/\Q)}c_{L}(E) )_{\pm}
\end{equation}
after the identification with an element of $H^{2}(K_{\beta}L/\Q,\{\pm 1\})$. The embedding problem associated with $\xi_{K_{\beta}L}$ is unobstructed and admits solutions $\gamma \in K_{\beta}$. Applying theorem 3.1 in \cite{quer2} will held one of these solutions, but first we need to restate our embedding problem in a way compatible with the notation of the theorem. Note that $c_{L}($E$)_{\pm} = c_{A,B}$ with $A=\{-10\}$ and $B=\{5\}$ where $c_{L}($E$)_{\pm}$ is given by the signs in table \ref{cociclo}. Let $\xi_{A,B}$ and $\theta_{\epsilon}$  be as in \cite{quer2}, where the character $\epsilon$ is the one above. Now using to the notation of theorem 3.1 in \cite{quer2} we consider the decomposition
$$K=\Q(\theta,\sqrt{-2}) = LMN=\Q(\theta)\Q(\sqrt{-10})\Q,$$
which gives $e=1$ ($L\cap M = \Q$), $m=1$ and $s=0$, and let $\sigma_{0},\sigma_{1} \in Gal(K/\Q)$ reflect this decomposition, that is
$$\left\{ \begin{array}{lll}
                     \sigma_{0}(\sqrt{5}) = -\sqrt{5}  & \mbox{ and } & \sigma_{0}(\sqrt{-2}) = -\sqrt{-2} \\
                     \sigma_{0}(\sqrt{5}) = \sqrt{5}  & \mbox{ and } & \sigma_{1}(\sqrt{-2}) = -\sqrt{-2}\
                  \end{array} \right..$$

We put $c=\xi_{A,B}\theta_{\epsilon}$ and compute $\zeta = c(1,\sigma_0)c(\sigma_0,\sigma_0)c(\sigma_0^2,\sigma_0)c(\sigma_0^{3},\sigma_0)=-1$. Theorem 3.1 in \cite{quer2} says that the embedding problem $(K/\Q, \{\pm 1\},[c])$ is solvable if and only if we can find elements $\alpha_{0},\alpha_{1}$ such that
$$N_{\sigma_{0}}(\alpha_0)=-1, \quad \quad N_{\sigma_{1}}(\alpha_1)=5, \quad \quad \mbox{ and } \frac{{^{\sigma_1}}\alpha_0}{\alpha_0} = \frac{{^{\sigma_0}}\alpha_1}{\alpha_1}.$$ Since $c = \xi_{K_{\beta}L}$ (here we are using notation of equation (\ref{soluvel})) we already know that this problem is unobstructed and a small search leads us to
$$ \alpha_0 = \frac{1}{2}(-1 + \theta + \theta^2)\sqrt{-2} \quad \quad \mbox{and} \quad \quad \alpha_1 = -5 + 2\theta^2.$$ Moreover, the same theorem gives us a splitting map for the cocycle $c$ given by
\begin{eqnarray*}
\beta_{\mbox{id}}&=&\beta_{\sigma_1}=1,\\
\beta_{\sigma_0} &=&\beta_{\sigma_1\sigma_0} = \alpha_0, \\
\beta_{\sigma_{0}^{2}} &=& \beta_{\sigma_1\sigma_{0}^{2}} = 2-\theta - \theta^2, \\
\beta_{\sigma_{0}^{3}} &=& \beta_{\sigma_1\sigma_{0}^{3}} = \frac{\sqrt{-2}}{2}(-\theta^2 - \theta + 3).\
\end{eqnarray*}
Finally, taking $x=1/4$ in formula (1) in \cite{quer2} we conclude that $\gamma = 2\theta^2 - \theta -5 \in K$ is a solution. Now a direct application of theorem 5.4 in \cite{quer1} to $E_{\gamma}$ yields theorem \ref{twist}. \qed

\bibliography{bibArtigo_1}

\begin{thebibliography}{10}

\bibitem{bil1}
N.~Billerey.
\newblock {\'E}quations de {F}ermat de {T}ype (5,5,p).
\newblock {\em Bull. Austral. Math. Soc.}, 76(2):161--194, 2007.

\bibitem{BD}
N.~Billerey and L.~Dieulefait.
\newblock Solving {F}ermat-type equations $x^{5} + y^{5} = d z^{p}$.
\newblock {\em Math. Comp.}, 79:535--544, 2010.

\bibitem{car}
H.~Carayol.
\newblock Sur les repr{\'e}sentation galoisiennes modulo $l$ attach{\'e}es aux
  formes modulaires.
\newblock {\em Duke Math. J.}, 59:785--801, 1989.

\bibitem{DJ}
L.~Dieulefait and J.~Jim{\'e}nez.
\newblock Solving {F}ermat-type equations $x^4 + d y^2 = z^p$ via modular
  $\mathbb{Q}$-curves over polyquadratic fields.
\newblock {\em J. Reine Angew. Math.}, 633:183--196, 2009.

\bibitem{ell2}
J.~Ellenberg.
\newblock Galois representations attached to $\mathbb{Q}$-curves and the
  generalized {F}ermat equation ${A}^4 + {B}^2 = {C}^p$.
\newblock {\em Amer. J. Math}, 126:763--787, 2004.

\bibitem{mil1}
J.~S. Milne.
\newblock On the arithmetic of abelian varieties.
\newblock {\em Invent. Math.}, 17:177--190, 1972.

\bibitem{pyle}
E.~E. Pyle.
\newblock Abelian varieties over $\mathbb{Q}$ with large endomorphisms algebras
  and their simple components over $\bar{\mathbb{Q}}$.
\newblock {\em Progress in Mathematics}, 224:189--234, 2004.

\bibitem{quer1}
J.~Quer.
\newblock $\mathbb{Q}$-curves and {A}belian varieties of
  $\mathit{GL}_{2}$-type.
\newblock {\em Proc. London Math. Soc.}, (3) 81:285--317, 2000.

\bibitem{quer2}
J.~Quer.
\newblock Embedding problems over abelian groups and an application to elliptic
  curves.
\newblock {\em J. Algebra}, 237:186--202, 2001.

\bibitem{rib1}
K.~Ribet.
\newblock Abelian {V}arieties over $\mathbb{Q}$ and modular forms.
\newblock {\em Progress in Mathematics}, 224:241--261, 2004.

\bibitem{sil2}
J.~H. Silverman.
\newblock {\em Advances {T}opics in the {A}rithmetic of {E}lliptic {C}urves}.
\newblock Springer-Verlag, 1994.

\end{thebibliography}

\bibliographystyle{plain}

\end{document}